\documentclass[12pt]{amsart}
\usepackage[latin1]{inputenc}
\usepackage{amsxtra,amssymb,amsthm,amsmath,amscd,mathrsfs, epsfig,eufrak,url}
\usepackage[all]{xy}
\def \C {{\mathbb C}}
\def \D {{\mathbb D}}

\def \N {{\mathbb N}}

\def \R {{\mathbb R}}

\def \Z {{\mathbb Z}}

\def \d {\,{\rm d}}
\def\re{{\Re e\,}}

\def\e{{\varepsilon}}

\def\le{\leqslant}

\def\ge{\geqslant}

\theoremstyle{plain}
\newtheorem{theorem}{Theorem}

\newtheorem{corollary}{Corollary}

\theoremstyle{remark}

\theoremstyle{definition}

\numberwithin{equation}{section}

\begin{document}

\vskip 5mm

\title[On the general additive divisor problem]
{On the general additive divisor problem}
\author{Aleksandar Ivi\'c  and Jie Wu}

\address{Katedra Matematike
\\
Rudarsko-geolo\v ski Fakultet
\\
Universitet u Beogradu
\\
 \DJ u-\break\v sina 7
\\
11000 Beograd
\\
Serbia}
\email{ivic@rgf.bg.ac.rs}

\address{
School of Mathematics
\\
Shandong University
\\
Jinan
\\
Shandong 250100
\\
China
\\
Institut Elie Cartan Nancy
\\
CNRS, Universit\'e Henri Poincar\'e (Nancy 1), INRIA
\\
Boulevard des Aiguillettes, B.P. 239
\\
54506 Van\-d\oe uvre-l\`es-Nancy
\\
France}
\email{wujie@iecn.u-nancy.fr}

\date{\today}

\begin{abstract}
We obtain a new upper bound for $\sum_{h\le H}\Delta_k(N,h)$
for $1\le H\le N$, $k\in \N$, $k\ge3$,
where $\Delta_k(N,h)$ is the (expected) error term in the asymptotic formula
for $\sum_{N<n\le2N}d_k(n)d_k(n+h)$,
and $d_k(n)$ is the divisor function generated by $\zeta(s)^k$.
When $k=3$ the result improves, for $H\ge N^{1/2}$, the bound given in the recent work \cite{[1]}
of Baier, Browning, Marasingha and Zhao,
who dealt with the case $k=3$.
\end{abstract}

\maketitle

\addtocounter{footnote}{1}

\centerline{\it 
In honor of Professor A.A. Karatsuba's 75th birthday}

\section{Introduction}

Let $d_k(n)$ denote that (generalized) divisor function, which represents
the number of ways $n$ can be written as a product of $k\;(\in\N)$ factors.
Thus
$$
\sum_{n=1}^\infty \frac{d_k(n)}{n^s} = \zeta(s)^k\qquad(\re s >1),
$$
where $\zeta(s)$ is the familiar zeta-function of Riemann.
In particular $d_1(n)\equiv 1$ and $d_2(n) = \sum_{\delta|n}1$
is the number of positive divisors of $n$.
The function $d_k(n)$ is a multiplicative function of $n$, and
$$
d_k(p^\nu) = (-1)^\nu{\binom{-k}{\nu}} = \frac{k(k+1)\cdots (k+\nu-1)}{\nu!}
$$
for primes $p$ and $\nu\in\N$.
The general divisor problem deals with the estimation of $\Delta_k(x)$,
the error term in the asymptotic formula
(see Chapter 13 of Ivi\'c \cite{[3]} and Chapter 12 of Titchmarsh \cite{[15]} for an extensive discussion)
\begin{equation}\label{(1.1)}
D_k(x) := \sum_{n\le x}d_k(n) = xp_{k-1}(\log x) + \Delta_k(x),
\end{equation}
where
\begin{equation}\label{(1.2)}
p_{k-1}(\log x)
= \mathop{{\rm Res}}_{s=1}\bigg(\zeta(s)^k \frac{x^{s-1}}{s}\bigg).
\end{equation}
Since $\zeta(s)$ is regular in $\C$ except at $s=1$ where it has a simple
pole with residue 1, it transpires
that $p_{k-1}(y)$ is a polynomial of degree $k-1$, whose coefficients
may be explicitly evaluated, and in particular $p_1(y) = y + 2\gamma -1$, where
$\gamma$ is Euler's constant. The important constants $\alpha_k, \beta_k$ are defined as
\begin{equation}\label{(1.3)}
\begin{aligned}
\alpha_k
& := \inf\big\{\,a_k\;:\;\Delta_k(x) \ll x^{a_k}\big\},
\\
\beta_k
& := \inf\bigg\{\,b_k\;:\;\int_1^X |\Delta_k(x)|^2\d x \ll  X^{1+2b_k}\bigg\}.
\end{aligned}
\end{equation}
It is known that $\alpha_k\ge \beta_k\ge (k-1)/(2k)$ for all $k\in \N$,
and the conjecture that
$\alpha_k= \beta_k= (k-1)/(2k)$  for all $k\in \N$ is equivalent to the Lindel\"of hypothesis
that $\zeta(\frac{1}{2}+{\rm i}t)\ll_\varepsilon (|t|+1)^\e$.
Here and later $\e \;(>0)$ denotes arbitrarily small constants, not necessarily
the same ones at each occurrence, while $\ll_{a,b,\ldots}$ means that the implied
constant in the $\ll$--symbol depends on $a,b,\ldots\,$.

\medskip

The general additive divisor problem is another important problem involving the divisor function $d_k(n)$.
It consists of the estimation of the quantity $\Delta_k(x,h)$, given by the formula
\begin{equation}\label{(1.4)}
\sum_{n \le x} d_k(n) d_k(n+h) = x\,P_{2k-2}(\log x;h)+\Delta_k(x,h).
\end{equation}
In \eqref{(1.4)} it is assumed that $k\ge2$ is a fixed integer,
and $P_{2k-2}(\log x;h)$ is a suitable polynomial of degree $2k-2$ in $\log x$,
whose coefficients depend on $k$ and $h$,
while $\Delta_k(x,h)$ is supposed to be the error term.
This means that we should have
\begin{equation}\label{(1.5)}
\Delta_k(x,h)\; =\; o(x) \quad{\rm as}\quad x \to \infty,
\end{equation}
but unfortunately \eqref{(1.5)} is not yet known to hold for any $k \ge 3$, even for fixed $h$.
However, when we consider the sum
\begin{equation}\label{(1.6)}
\sum_{h \le H}\Delta_k(x,h),
\end{equation}
we may reasonably hope that a certain cancelation will occur
among the individual summands $\Delta_k(x,h)$,
since there are no absolute value signs in \eqref{(1.6)}.
It turns out that it is precisely the estimation of the sum in \eqref{(1.6)}
which is relevant for bounding the integral
\begin{equation}\label{(1.7)}
\int_0^T|\zeta(\tfrac{1}{2}+{\rm i}t)|^{2k}\d t,
\end{equation}
which is of great importance in the theory of the Riemann zeta-function (see
the monographs \cite{[3], [4], [15]}).

\medskip

For $k=1$ the sum in \eqref{(1.6)} is trivial,
while for $k=2$ it was extensively studied by many authors,
including Kuznetsov \cite{[10]}, Motohashi \cite{[13]},
Ivi\'c \& Motohashi \cite{[8]} and Meurman \cite{[12]}.
The natural next step in \eqref{(1.6)} is  to deal with the case $k=3$,
but the works of A.I.~Vinogradov and Takhtad\v zjan \cite{[19], [20]} and A.I.~Vinogradov \cite{[16], [17], [18]}
show that the analytic problems connected with the Dirichlet series generated
by $d_3(n)d_3(n+h)$ are overwhelmingly hard.
The ensuing problems are connected with the group $SL(3, \Z)$,
and they are much more difficult than the corresponding problems connected with the group
$SL(2, \Z)$ which appear in the case $k=2$.
The latter involve the spectral theory of the non-Euclidean Laplacian,
which was extensively developed in recent times by Kuznetsov (see e.g., \cite{[11]}),
Iwaniec and others (see Motohashi's monograph \cite{[14]} for applications of spectral
theory to the theory of $\zeta(s)$).
Thus at present in the case $k=2$ we have sharp explicit formulas, while in the case $k>2$
we have none.

\medskip

A.I.~Vinogradov \cite{[18]} conjectured that $\Delta_k(x,h) \ll x^{1-1/k}$,
without stating for which range of $h$ this sharp bound should hold.
Very likely this bound is too strong, and  (even for fixed $h$)
it seems probable that a power of a logarithm should be included on the right-hand side.
More importantly, one hopes that the bound
\begin{equation}\label{(1.8)}
\sum_{h\le H}\Delta_k(x,h) \ll_{k,\e} Hx^{1-1/k+\e}
\quad{{\rm for}}\quad 1 \le H \le x^{(k-2)/k +\delta_k}
\end{equation}
holds uniformly in $H$ for fixed $k \ge 3$ and some $\delta_k > 0$,
which was stated in \cite{[5]}.
Note that Vinogradov's conjecture in the form
$\Delta_k(x,h) \ll_{k,\e} x^{1-1/k+\e}$ trivially implies \eqref{(1.8)}, but the
important point is that there are no absolute value signs in the sum
in \eqref{(1.8)}.  One can also assume \eqref{(1.8)} to hold
in the case $k=2$ for $1 \le H \le \sqrt x$, say.
Then it would follow that the inequality
\begin{equation}\label{(1.9)}
\int_{T-G}^{T+G} \Bigl|\zeta(\tfrac{1}{2}+{\rm i}t)\bigr|^4 \d t
\ll_\varepsilon GT^\e
\end{equation}
holds with $G=T^{5/6}$, whereas it is known (see e.g., \cite{[7]}) that $G=T^{2/3}$
is unconditionally permissible. It was conjectured in \cite{[5]} that for any
$k\ge 2$ and $h \ge 1$ one has
\begin{equation}\label{(1.10)}
\Delta_k(x,h) \, =\, \Omega\bigl(x^{1-1/k}\bigr).
\end{equation}
For $k=2$ and fixed $h$ this conjecture was proved by Motohashi \cite{[13]}.
As usual, $f(x) = \Omega\bigl(g(x)\bigr)$ means that
$\lim\limits_{x\to\infty}\, {f(x)/g(x)} \ne 0$.

\smallskip

The general additive divisor problems is connected to the power moments
of $|\zeta(\frac{1}{2}+{\rm i}t)|$ (see e.g., \cite{[3]} and \cite{[4]} for an extensive account).
In 1996  the first author \cite{[5]} proved that
\begin{equation}\label{(1.11)}
\int_0^T \bigl| \zeta(\tfrac{1}{2}+{\rm i}t) \bigr|^6 \d t
\ll_\varepsilon T^{1+\e}+ T^{(\alpha+3\beta-1)/2+\e}
\end{equation}
provided that
$$
\sum_{h\le H}\Delta_3(x,h) \ll_\varepsilon H^\alpha x^{\beta+\varepsilon}
$$
holds for $1\le H\le x^{1/3+\delta_3}$ for some constant $\delta_3>0,
0\le\alpha, \, \beta\le1, \alpha+\beta\ge1$.
The conjecture \eqref{(1.8)} with $k=3$ means that we can take $\alpha=1,\beta=2/3$ in \eqref{(1.11)}
so that the sixth moment in the form
\begin{equation}\label{(1.12)}
\int_0^T \bigl|\zeta(\tfrac{1}{2}+{\rm i}t)\bigr|^6 \d t
\ll_\varepsilon\; T^{1+\e}
\end{equation}
follows.
Note that the best known exponent
of $T$ for the right-hand side of the above integral is 5/4 (see \cite[Chapter 8]{[3]}).

\medskip

In \cite{[6]} the research begun in \cite{[5]} was continued, and a plausible heuristic
evaluation of the polynomial $P_{2k-2}(x;h)$ in \eqref{(1.4)} was made.
Yet another (heuristic) evaluation of the sum in \eqref{(1.5)}
was made  later by Conrey and Gonek \cite{[2]} in 2001. Moreover, it was
shown in \cite{[6]} that, for a fixed integer $k \ge 3$ and any fixed $\e > 0$, we have
\begin{equation}\label{(1.13)}
\int_0^T|\zeta(\tfrac{1}{2}+{\rm i}t)|^{2k}\d t
\ll_{k,\e} T^{1+\e}\bigg(1 + \sup_{T^{1+\e}<M\ll T^{k/2}} \frac{G_k(M;T)}{M}\bigg),
\end{equation}
if, for $T^{1+\e} \le M \ll T^{k/2}$ and $M < M' \le 2M$,
\begin{equation}\label{(1.14)}
G_k(M;T)
:= \sup_{\substack{M\le x \le M'\\ 1\le t \le M^{1+\e}/T}} \Big|\sum_{h\le t}\,\D _k(x,h)\Big|.
\end{equation}
This result, which generalizes \eqref{(1.11)},
provides a directlink between upper bounds for the $2k$-thmoment of $|\zeta(\frac{1}{2}+{\rm i}t)|$
and sums of $\D _k(x,h)$ over the shift parameter $h$.
The result also gives an insight as to the limitations of the attack
on the $2k$-th moment of $|\zeta(\frac{1}{2}+{\rm i}t)|$ via the use of estimates for $\Delta_k(x,h)$.
Of course the problem greatly increases in complexity as $k$ increases,
and this is one of the reasons why in \cite{[5]} only the case $k =3$ was
considered. The case $k = 2$ was not treated, since for the fourth
moment of $|\zeta(\frac{1}{2}+{\rm i}t)|$ we have an asymptotic formula with precise results for
the corresponding error term (see e.g.,  \cite{[7]} and \cite{[14]}).
Note that \eqref{(1.13)}--\eqref{(1.14)} again lead to the sixth moment
bound \eqref{(1.12)} if the conjecture \eqref{(1.8)} holds with $k=3$.

\vskip 8mm

\section{The general additive divisor problem}

The main objective of this note is to study  the averaged  sum \eqref{(1.6)},
when $k\ge3$ is a fixed integer.
To this end we introduce more notation, defining
\begin{equation}\label{(2.1)}
D_k(N,h) := \sum_{N<n\le2N}d_k(n)d_k(n+h),
\end{equation}
and letting henceforth
\begin{equation}\label{(2.2)}
\Delta_k(N;h) := D_k(N,h) - \int_N^{2N}{\mathfrak S}_k(x,h)\d x,
\end{equation}
so that $\Delta_k(N;h)$ in \eqref{(2.2)} differs slightly from \eqref{(1.4)};
in fact it equals $\Delta_k(2N,h) - \Delta_k(N,h)$ in the notation of \eqref{(1.4)}.
Here we follow the notation of \cite{[1]}, based on the approach of
Conrey and Gonek \cite{[2]}, who made conjectures on the high moments
of $|\zeta(\tfrac{1}{2}+{\rm i}t)|$. Let us also define
$$
{\mathfrak S}_k(x,h) := \sum_{q=1}^{\infty} \frac{c_q(h)}{q^2} Q_k(x,q)^2,
$$
where $\mu(n)$ is the M\"obius function,
$c_q(h):=\sum_{d\mid (h,q)} d\mu(q/d)$ is the Ramanujan sum and
$Q_k(x,q)$ is defined as follows.
If $\varphi(n)$ is the Euler totient function, set
$$
\Psi_{d, e}(s, q, k)
:= \frac{d\mu(d)\mu(e)}{\varphi(d)e}
\prod\limits_{p\mid (eq/d)} \bigg\{\bigg(1-\frac{1}{p^s}\bigg)^k
\sum_{\nu=0}^{\infty} \frac{d_k\left(p^{\nu+\nu_p(eq/d)}\right)}{p^{\nu s}}\bigg\},
$$
where here and later $\nu_p(m)$ is the $p$-adic valuation of $m$. Then we define
\begin{align*}
Q_k(x,q)
& :=
\frac{1}{2\pi {\rm i}} \int\limits_{|s-1|=1/8} \zeta(s)^k
\sum_{d\mid q} \sum_{e\mid d} \Psi_{d, e}(s, q, k)
\bigg(\frac{ex}{dq}\bigg)^{s-1} \d s
\\
& \; =
\mathop{{\rm Res}}\limits_{s=1} \,
\bigg\{\zeta(s)^k \sum_{d\mid q} \sum_{e\mid d} \Psi_{d, e}(s, q, k)
\bigg(\frac{ex}{dq}\bigg)^{s-1}\bigg\},
\end{align*}
by the residue theorem. Thus $Q_k(x,q)$ is a polynomial of degree $2k-2$
whose coefficients depend on $q$, and may be explicitly evaluated.
The work of Conrey-Gonek (op. cit.) predicts, as stated in (2.2), that
$D_k(N,h)$ is well approximated by $\int_N^{2N}{\mathfrak S}_k(x,h)\d x$, which equals
$N$ times a polynomial in $\log N$ of degree $2k-2$,
all of whose coefficients depend on $h$ and $k$.
This is in agreement with  [5] (when $k=3$ and [6] (in the general case), although the shape of the
polynomial in question is somewhat different. Conrey and Gonek even predict that
uniformly
$$
\Delta_k(N; h)
\ll_\varepsilon N^{1/2+\e}\quad
{{\rm for}}\quad 1 \le h \le N^{1/2}.
$$
This conjecture is probably too strong,
and one feels that more likely the bound $\Delta_k(N; h) \ll_\varepsilon N^{1-1/k+\e}$
is closer to the truth (see (1.10)).

\medskip

In a recent work \cite{[1]}, Baier, Browning, Marasingha and Zhao obtain
new results involving averages of $\Delta_3(N; h)$ (they employ the terminology
``shifted convolutions of $d_3(n)$'', which seems appropriate). They proved that
\begin{equation}\label{(2.3)}
\sum_{h\le H} \Delta_3(N; h)
\ll_\varepsilon
N^{\varepsilon}\big(H^2+H^{1/2}N^{13/12}\big)
\qquad(1\le H\le N),
\end{equation}
and if $N^{1/3+\varepsilon}\le H\le N^{1-\varepsilon}$, then
there exists $\delta(\e)>0$ such that
\begin{equation}\label{(2.4)}
\sum\limits_{h\le H} \left| \Delta_3(N; h) \right|^2
\ll_\varepsilon HN^{2-\delta(\e)}.
\end{equation}
These results can be used, in conjunction with the bounds \eqref{(1.13)}--\eqref{(1.14)} when $k=3$,
to bound the integral in \eqref{(1.11)},
but they will produce only the exponent 11/8 on the right-hand side of \eqref{(1.11)},
hence no improvement on the known result for the sixth moment of $|\zeta(\tfrac12+{\rm i}t)|$.

\medskip

{\bf Remark 1}.
Note that \eqref{(2.3)}, in the range $N^{1/6+\e}\le H \le N^{1-\e}$, provides
an asymptotic formula for the averaged sum $\sum_{h\le H} D_3(N,h)$ (see \eqref{(2.1)}).
However, it should be noted that no asymptotic formula for the individual $D_3(N,h)$ has been
found yet, and in general for $\Delta_k(N; h)$ when $k\ge3$. In fact, it is worth
pointing out that when $1\le H\le N^{1/6}$, the bound in \eqref{(2.3)} is worse than
the trivial bound $HN^{1+\varepsilon}$. Namely we have
$$
\sum_{h\le H} D_k(N, h)
\ll_\varepsilon \sum_{h\le H} \sum_{N<n\le 2N} (n+h)^{\varepsilon/2}
\ll_\varepsilon (HN)^{1+\varepsilon/2}
\ll_\varepsilon HN^{1+\varepsilon}.
$$
On the other hand we have
$$
\sum_{h\le H} \int_{N}^{2N} {\mathfrak S}_k(x, h) \d x
\,\ll_\varepsilon HN^{1+\varepsilon},
$$
which is obvious from \eqref{(3.8)}.
Hence by \eqref{(2.2)} it follows that
\begin{equation}\label{(2.5)}
\sum_{h\le H}\Delta_k(N; h)
\ll_\varepsilon HN^{1+\e}\qquad(1\le H\le N),
\end{equation}
and clearly \eqref{(2.5)} for $k=3$ improves \eqref{(2.3)} for $1\le H\le N^{1/6}$.
\medskip
The aim of this note is to give a bound for the sum in \eqref{(1.6)}, or equivalently
for the average of \eqref{(2.2)},
which for $k=3$ improves \eqref{(2.3)} for a certain range of $H$.
The result is contained in the following

\begin{theorem}\label{thm1}
For fixed $k\ge 3$ we have
\begin{equation}\label{(2.6)}
\sum_{h\le H}\Delta_k(N; h)
\ll_\varepsilon N^\e \big(H^2 + N^{1+\beta_k}\big)
\qquad(1\le H\le N),
\end{equation}
where $\beta_k$ is defined by \eqref{(1.3)}.
\end{theorem}

Note that we have $\beta_3 = 1/3, \beta_4 = 3/8$ (see Chapter 13 of \cite{[3]}),
$\beta_5 \le 9/20$ (see Zhang \cite{[21]}), $\beta_6 \le 1/2$, etc.
For a discussion of the values of $\alpha_k$
and $\beta_k$, see also the paper by Ouellet and Ivi\'c \cite{[9]}.

\begin{corollary}.
We have, for $1 \le H\le N$,
\begin{equation}\label{(2.7)}
\begin{aligned}
\sum_{h\le H}\Delta_3(N; h) &\ll_\varepsilon N^\e \big(H^2 + N^{4/3}\big),
\\
\sum_{h\le H}\Delta_4(N; h) &\ll_\varepsilon N^\e \big(H^2 + N^{11/8}\big),
\\
\sum_{h\le H}\Delta_5(N; h) &\ll_\varepsilon N^\e \big(H^2 + N^{29/20}\big),
\\
\sum_{h\le H}\Delta_6(N; h) &\ll_\varepsilon N^\e \big(H^2 + N^{3/2}\big).
\end{aligned}
\end{equation}
\end{corollary}

\medskip

{\bf Remark 2}.
Since it is known that $\beta_k <1$ for any $k$, this means
that the bound in \eqref{(2.6)} improves on the trivial bound $HN^{1+\e}$ in the
range $N^{\beta_k+\e} \le H\le N^{1-\e}$. Our result thus  supports
the assertion that $\Delta_k(N; h)$ is really the error term in the asymptotic
formula for $D_k(N,h)$, as given by \eqref{(3.1)} and \eqref{(3.2)}.
In the case when $k=3$, we have
by \eqref{(2.7)} an improvement of \eqref{(2.3)} when $H\ge N^{1/2}$.

\vskip 8mm

\section{Proof of Theorem \ref{thm1}}

We begin by noting that obviously
$$
\sum_{h\le H}d_k(n+h) = \sum_{m\le n+H}d_k(m) - \sum_{m\le n}d_k(m).
$$
Therefore by \eqref{(1.1)}--\eqref{(1.2)} and \eqref{(2.1)}--\eqref{(2.2)}  we can write
\begin{equation}\label{(3.1)}
\begin{aligned}
\sum_{h\le H} \Delta_k(N, h)
& = \sum_{N<n\le 2N} d_k(n) \sum_{h\le H} d_k(n+h) - \sum_{h\le H} \int_N^{2N} {\mathfrak S}_k(x, h) \d x
\\
& = M_k(N, H) + R_k(N, H) - \sum_{h\le H} \int_N^{2N} {\mathfrak S}_k(x, h) \d x,
\end{aligned}
\end{equation}
say, where
\begin{align*}
M_k(N, H)
& := \sum_{N<n\le 2N} d_k(n) \mathop{{\rm Res}}_{s=1}\bigg(\zeta(s)^k \frac{(n+H)^s-n^s}{s}\bigg),
\\
R(N, H)
& := \sum_{N<n\le 2N} d_k(n) \big(\Delta_k(n+H)-\Delta_k(n)\big),
\end{align*}
where $\Delta_k(x)$ is defined by \eqref{(1.1)}.
It is rather easy to estimate $R_k(N, H)$. Namely
since $d_k(n)\ll_{\varepsilon} n^{\varepsilon}$, we have trivially
$$
R_k(N,H) \ll_\varepsilon N^{\varepsilon} \sum_{n\le 3N} |\Delta_k(n)|.
$$
For $n<t<n+1$, we obviously have
$$
\Delta_k(n) - \Delta_k(t)
= tp_{k-1}(\log t)-np_{k-1}(\log n)
\ll (\log n)^{k-1}.
$$
Thus
\begin{equation}\label{(3.2)}
\begin{aligned}
R_k(N, H)
& \ll_\varepsilon N^{\varepsilon} \sum_{n\le 3N} \int_n^{n+1} |\Delta_k(n)| \d t
\\
& \ll_\varepsilon N^{\varepsilon} \sum_{n\le 3N} \int_n^{n+1} |\Delta_k(t)| \d t + N^{1+\varepsilon}
\\
& \ll_\varepsilon N^{\varepsilon} \int_1^{4N} |\Delta_k(t)| \d t + N^{1+\varepsilon}
\\
& \ll_\varepsilon N^{\varepsilon} \bigg(N \int_1^{4N} |\Delta_k(t)|^2 \d t\bigg)^{1/2} + N^{1+\varepsilon}
\\\noalign{\vskip 1,5mm}
& \ll_\varepsilon N^{1+\beta_k+\varepsilon},
\end{aligned}
\end{equation}
where we used the Cauchy-Schwarz inequality for integrals and the mean square bound
\eqref{(1.3)} in the last step.

\medskip
To estimate $M_k(N,H)$, set
$$
u_k(x) := \mathop{{\rm Res}}_{s=1}\bigg(\zeta(s)^k \frac{(x+H)^s-x^s}{s}\bigg).
$$
Then we can write
$$
M_k(N,H) = \int_N^{2N+0}u_k(x)\d D_k(x).
$$
But we have, since
$$
D_k(x) = \mathop{{\rm Res}}_{s=1}\bigg(\zeta(s)^k \frac{x^s}{s}\bigg)
+\Delta_k(x)
$$
in view of \eqref{(1.1)} and \eqref{(1.2)},
\begin{equation}\label{(3.3)}
M_k(N,H) = \int_N^{2N}u_k(x)\mathop{{\rm Res}}_{s=1}\big(\zeta(s)^k x^{s-1}\big)\d x
+ \int_N^{2N}u_k(x)\d\Delta_k(x).
\end{equation}
Further note that
\begin{equation}\label{(3.4)}
\begin{aligned}
u_k(x)
& = yp_{k-1}(\log y)\Bigl|_x^{x+H} \ll H(\log x)^{k-1},
\\
u'_k(x)
& = \mathop{{\rm Res}}_{s=1}\zeta(s)^k
\big\{(x+H)^{s-1}-x^{s-1}\big\}
\ll_\varepsilon x^\e.
\end{aligned}
\end{equation}
On integrating by parts and using \eqref{(1.3)} and \eqref{(3.4)} we obtain, similarly to \eqref{(3.2)},
\begin{equation}\label{(3.5)}
\begin{aligned}
\int_N^{2N}u_k(x)\d\Delta_k(x)& = u_k(x)\Delta_k(x)\Bigl|_N^{2N} - \int_N^{2N}\Delta_k(x)u'_k(x)\d x
\\
& \ll_\varepsilon HN^{\alpha_k+\e} + N^{1+\beta_k+\e}.
\end{aligned}
\end{equation}
As for the other integral in \eqref{(3.3)}, note that
$$
\frac{(x+H)^s-x^s}{s}
= \frac{x^s}{s} \bigg\{1+\frac{sH}{x} + \binom{s}{2} \frac{H^2}{x^2}+\cdots-1\bigg\}.
$$
This gives
\begin{equation}\label{(3.6)}
\int_N^{2N}u_k(x)\mathop{{\rm Res}}_{s=1}\big(\zeta(s)^k x^{s-1}\big)\d x
= H\int_N^{2N}\big(\mathop{{\rm Res}}_{s=1}\zeta(s)^kx^{s-1}\big)^2\d x
+ O_\varepsilon\big(H^2N^\e\big).
\end{equation}
Therefore from \eqref{(3.3)}, \eqref{(3.5)} and \eqref{(3.6)} we obtain
\begin{equation}\label{(3.7)}
\begin{aligned}
M_k(N,H)
& = H\int_N^{2N}\big(\mathop{{\rm Res}}_{s=1}
\zeta(s)^kx^{s-1}\big)^2\d x
\\
& \quad
+ O_\varepsilon\big(H^2N^\e + NH^{\alpha_k+\e} + N^{1+\beta_k+\e}\big).
\end{aligned}
\end{equation}

Next we shall prove that
\begin{equation}\label{(3.8)}
\sum\limits_{h\le H} \int_{N}^{2N} \mathfrak{S}_k(x,h)\d x
= H \int_{N}^{2N} \big(\mathop{{\rm Res}}
\limits_{s=1} \zeta(s)^kx^{s-1}\big)^2 \d x
+O_\varepsilon\big(N^{1+\varepsilon}\big).
\end{equation}
The case of $k=3$ has been treated in \cite{[1]}.
Here
we repeat the same argument with some simplification in the general case,
obtaining \eqref{(3.8)}.

First write
$$
x^{s-1}
= \sum_{n=0}^\infty \frac{(\log x)^n}{n!} (s-1)^n.
$$
Since $\Psi_{d, e}(s, q)$ and $(s-1)^n\zeta(s)^k$ with $n\ge k$ are holomorphic for $\re s>0$,
Cauchy's theorem allows us to deduce that
$$
Q_k(x,q)
= \frac{1}{2\pi {\rm i}} \sum_{n=0}^{k-1} \int\limits_{|s-1|=1/8} \zeta(s)^k
\sum_{d\mid q} \sum_{e\mid d}
\Psi_{d, e}(s, q) \frac{(\log(dx/eq))^n}{n!}(s-1)^n \d s.
$$
Clearly for $\re s>\tfrac{1}{2}$, we have
\begin{align*}
\Psi_{d, e}(s, q)
& \ll \frac{d}{\varphi(d)e}
\prod_{p\mid (eq/d)} \bigg\{\bigg(1+\frac{1}{p^{1/2}}\bigg)^k
\sum_{\nu=0}^\infty \frac{d_k(p^{\nu+\nu_p(eq/d)}}{p^{\nu s}}\bigg\}
\\
& \ll_\varepsilon \frac{d}{\varphi(d)e}
\prod_{p\mid (eq/d)} \bigg\{\bigg(1+\frac{1}{p^{1/2}}\bigg)^k
p^{\nu_p(eq/d)\e/4}\sum_{\nu\ge 0} \frac{p^{\nu\e/4}}{p^{\nu/2}}\bigg)
\bigg\}
\\\noalign{\vskip 1,5mm}
& \ll_\varepsilon q^{\e/2}.
\end{align*}
Thus
\begin{equation}\label{(3.9)}
Q_k(x,q)
\ll_{\e, k} q^{\e} (\log x)^{k-1},
\end{equation}
where the implied constant depends only on $\e$ and $k$.

\medskip
In view of (3.9) and the bound $|c_q(h)|\le (h, q)$,
we have
\begin{equation}\label{(3.10)}
\begin{aligned}
\sum_{h\le H} \sum_{q>H} \frac{c_q(h)}{q^2} Q_k(x, q)^2
& \ll (\log x)^{k-1} \sum_{h\le H} \sum_{q>H} \frac{(h, q)}{q^{2-\e}}
\\
& \ll_{\e, k} H^{\e} (\log x)^{k-1}.
\end{aligned}
\end{equation}
On the other hand, it is well known that
$
\sum_{h\le q} c_q(h) = 0
$
if $q>1$.
From this it is easy to deduce that
$$
\sum_{h\le H} c_q(h) = \begin{cases}
H + O(1)                              & \text{if $\;q=1$,}
\\
O_\varepsilon(q^{1+\e})   & \text{if $\;q>1$.}
\end{cases}
$$
With the help of this relation and (3.9), we can write
\begin{equation}\label{(3.11)}
\begin{aligned}
& \sum_{h\le H} \sum_{q\le H} \frac{c_q(h)}{q^2} Q_k(x, q)^2
\\
& = \{H + O(1)\} Q_k(x, 1)^2
+ O\bigg((\log x)^{k-1} \sum_{1<q\le H} \frac{1}{q^{1-\e}}\bigg)
\\
& = H\big(\mathop{{\rm Res}}\limits_{s=1} \zeta(s)^kx^{s-1}\big)^2
+ O\big((\log x)^{k-1} H^{\e}\big),
\end{aligned}
\end{equation}
where we have used the fact that
$$
Q_k(x, 1)
= \mathop{{\rm Res}}\limits_{s=1}\big(\zeta(s)^kx^{s-1}\big)
\ll_k (\log x)^{k-1}.
$$
By combining (3.10) and (3.11), we obtain (3.8).

\medskip
From (3.1), (3.2), (3.7) and (3.8) we obtain
\begin{equation}\label{(3.12)}
\sum_{h\le H}\Delta_k(N,H)
\ll_\varepsilon N^\e \big(H^2+ HN^{\alpha_k} + N^{1+\beta_k}\big)
\quad
(1\le H\le N).
\end{equation}
But we always have
\begin{equation}\label{(3.13)}
\alpha_k\le \tfrac{1}{2} + \tfrac{1}{2}\beta_k.
\end{equation}
To see this note that, for $1\le H\le x$, the defining relation (1.1) and $d_k(n) \ll_\varepsilon n^\e$
give
\begin{align*}
\Delta_k(x) - \frac{1}{H}\int_x^{x+H}\Delta_k(y)\d y
& = \frac{1}{H}\int_x^{x+H}\left(\Delta_k(x)-\Delta_k(y)\right)\d y
\\
&
\ll_\varepsilon \frac{1}{H} \int_x^{x+H}\big\{|D_k(x)-D_k(y)|+O(x^\e)\big\}\d y
\\\noalign{\vskip 1mm}
& \ll_\varepsilon Hx^\varepsilon.
\end{align*}
This gives, by the Cauchy-Schwarz inequality for integrals and \eqref{(1.3)},
\begin{align*}
\Delta_k(x)
& \ll_\varepsilon \frac{1}{H}\int_x^{x+H}|\Delta_k(y)|\d y + Hx^\e
\cr
& \ll_\varepsilon x^{1+\beta_k+\e}H^{-1} + Hx^\e
\\\noalign{\vskip 1mm}
& \ll_\varepsilon x^{(1+\beta_k)/2+\e}
\end{align*}
with $H = x^{(1+\beta_k)/2}$. Hence
$$
\Delta_k(x) \;\ll_\varepsilon\;x^{(1+\beta_k)/2+\e}
$$
and \eqref{(3.13)} follows.
Now in \eqref{(3.12)} we have $HN^{\alpha_k} \le H^2$ for $H\ge N^{\alpha_k}$.
If $H\le N^{\alpha_k}$, then $HN^{\alpha_k} \le N^{2\alpha_k} \le N^{1+\beta_k}$ by \eqref{(3.13)}.
Thus the term $HN^{\alpha_k}$ in \eqref{(3.12)} can be discarded, and \eqref{(2.6)} follows.
This completes the proof of the Theorem.

\end{document}